# Scalable Method for Linear Optimization of Industrial Processes

Leonid B. Sokolinsky, Irina M. Sokolinskaya
South Ural State University
Chelyabinsk, Russian Federation
leonid.sokolinsky@susu.ru, irina.sokolinskaya@susu.ru

*Abstract* — In the development of industrial digital twins, the optimization problem of technological and business processes often arises. In many cases, this problem can be reduced to a large-scale linear programming (LP) problem. The article is devoted to the new method for solving large-scale LP problems. This method is called the "apex-method". The apex-method uses the predictor-corrector framework. The predictor step calculates a point belonging to the feasible region of LP problem. The corrector step calculates a sequence of points converging to the exact solution of the LP problem. The article gives a formal description of the apex-method and provides information about its parallel implementation in C++ language by using the MPI library. The results of large-scale computational experiments on a cluster computing system to study the scalability of the apex method are presented.

## I. INTRODUCTION

In the development of industrial digital twins, the optimization problem of technological and business processes often arises. In many cases, this problem can be reduced to a large-scale linear programming (LP) problem. Known software often does not allow solving such large-scale LP problems in an acceptable time [1]. At the same time, the exascale computing systems must be established in the nearest 2-3 years [2], which are potentially capable to solve such huge problems. Therefore, the issue of developing new efficient methods for solving large-scale LP problems by exascale computing systems is important. Until now, one of the most common ways to solve the LP problem was a class of algorithms proposed and developed by Danzig based on the simplex method [3]. The simplex method was found to be efficient for solving a large class of LP problems. In particular, the simplex method easily takes advantage of any hyper-sparsity in a LP problem [4] However, the simplex method has some fundamental features that limit its use for the solution of large LP problems. First, in certain cases, the simplex method has to iterate through all the vertices of the simplex, which corresponds to exponential time complexity [5]–[7]. Second, in most cases, the simplex method successfully solves LP problems containing up to 50,000 variables. However, a loss of precision is observed when the simplex method is used for solving large LP problems [8]. Such a loss of precision cannot be compensated even by applying such computational intensive procedures as "affine scaling" or "iterative refinement" [9]. Third, in the general case, the sequential nature of simplex method makes it difficult to parallelize on multiprocessor systems with distributed memory [10]. Numerous attempts have been made to make a scalable parallel implementation of the simplex method, but all of them were unsuccessful [11]. In all cases, the scalability boundary was from 16 to 32 processor nodes (see, for example, [12]). Khachiyan proved [13], using a variant of the ellipsoid method (proposed in the 1970's by Shor [14] and Yudin and Nemirovskii [15]), that LP problems can be solved in polynomial time. However, attempts to apply this approach in practice were unsuccessful, since in the vast majority of cases, the ellipsoid method demonstrated much worse efficiency compared to the simplex method. Later, Karmarkar [16] proposed a polynomial-time interior-point algorithm, which can be used in practice. This algorithm generated the whole field of modern interior-point methods [17], which are able to solve large-scale LP problems with millions of variables and millions of equations [18]–[22]. Moreover, these methods are self-correcting, and therefore provide high precision of calculations. A general lack of interior-point methods is the need to find some feasible point satisfying all the constraints of the LP problem before starting calculations. Finding such an interior point can be reduced to solving an additional LP problem [23]. Another method for finding an interior point is the pseudo-projection method [24], which uses Fejer mappings [25]. One more significant disadvantage of the interior-point method is its poor scalability on multiprocessor systems with distributed memory. There are several successful parallel implementations of the interior-point method for particular cases (see, for example, [26]), but, in the general case, an efficient parallel implementation on exascale multiprocessor systems could not be built. In accordance with this, the research direction related to the development of new scalable methods for solving LP problems is urgent.

In the article [27], the authors proposed an idea of a new scalable iterative method based on the predictor-corrector scheme for solving large-scale LP problems on cluster computing systems. The method proceeds in two stages: Quest (predictor) and Target (corrector). The Quest stage searches for a feasible point that satisfies all the constraints of the LP problem. The Target stage, starting from the inner point obtained in the Quest stage, builds a sequence of feasible points that converges to the exact solution of the LP problem. The Quest stage was studied in [24], [27]–[29]. In [24], the notion of pseudo-projection on a convex closed set, generalizing the concept of projection, was introduced. The pseudo-projection method is used in the Quest stage to find an initial feasible point. To calculate the pseudo-projection, we use Fejer approximations [25], which can "self-correct" when the numerical error is accumulated. Paper [28] demonstrates that multi-core accelerators can be efficiently used for computing the pseudo-projection on high-dimensional polytope. In [27], the convergence theorem of the pseudo-projection algorithm was proved in the case when the polytope bounding the set of feasible points is shifted by a translation. In [29], a parallel algorithm implementing Target stage was proposed. This algorithm forms a special system of points having the shape of the n-dimensional axisymmetric cross. The cross moves in the n-dimensional space in such a way that the solution of the LP problem permanently was in the $\varepsilon$-vicinity of the central point of the cross. However, the mentioned algorithm has the time complexity exponentially depending on the problem dimension.

In this paper, we propose and investigate a new scalable iterative method for solving large-scale LP problems intended for cluster computing systems. We named this method the "apex method". The apex method also proceeds in two stages – Quest and Target, but it has the polynomial complexity unlike the previous approach. The paper is organized as follows. In section II, we give a mathematical formulation of the LP problem and the explanation of apex method. In section III, we briefly discuss a parallel implementation of the Target stage and present the results of large-scale computational experiments evaluating scalability of the apex method on a cluster computing system. Finally, in Section IV, we give our conclusions and comment on possible further investigations of the apex method.

## II. APEX METHOD

Consider the LP problem

$$\bar{x} = \arg\max_{x}\{\langle c,x\rangle \mid Ax \leq b\}, \quad (1)$$

where $\bar{x}, x, c \in \mathbb{R}^n$, $A \in \mathbb{R}^{n \times m}$ and $b \in \mathbb{R}^m$. Here, $\langle c, x \rangle$ denotes the dot product of $c$ and $x$. We assume that the constraint $x \geq 0$ is also included in the system $Ax \leq b$. Let us denote by $a_i$ the $i$-th row of the matrix $A$. We assume from now on that $a_i$ not equal to the zero vector for all $i = 1,\ldots,m$. Let $M$ denote the $n$-dimensional polytope bounding the set of feasible points of the problem (1). Such polytope is always a convex closed set. We also assume the polytope $M$ to be bounded. By definition, $x$ is a *boundary point* of the polytope $M$ if any vicinity of x has non-empty intersection with both, $M$ and its complement. Let us define $\Gamma_M$, the boundary of polytope $M$, as the set of all its boundary points.

The *apex method* uses predictor-corrector scheme and proceeds in two stages: *Quest* (predictor) and *Target* (corrector). The *Quest stage* searches for a feasible point $\tilde{x} \in M$. The *Target stage*, using $\tilde{x}$, calculates a sequence of points $\{u_0, u_1, \ldots, u_k, \ldots\}$ that has the following properties:

$$u_k \in \Gamma_M; \quad (2)$$

$$\langle c, u_k \rangle < \langle c, u_{k+1} \rangle; \quad (3)$$

$$\lim_{k \to \infty} \langle c, u_k \rangle = \bar{x}. \quad (4)$$

The condition (2) means that all points in the sequence "lie" on the boundary of the polytope $M$. The condition (3) means that the value of the objective function at each next point is greater than at the previous one. The condition (4) means that the sequence converges to the exact solution of the problem (1). The Quest stage finds the point $\tilde{x} \in M$ using the pseudo-projection operation [24], which is a generalization of the orthogonal projection operation. Let us give a formal definition. Let the mapping $\varphi_M : \mathbb{R}^n \to \mathbb{R}^n$ be given by

$$\varphi_M(x) = \frac{1}{h}\sum_{i=1}^{m}\rho_i^+(x), \quad (5)$$

where

$$\rho_i^+(x) = \frac{\max\{\langle a_i, x\rangle - b_i, 0\}}{\|a_i\|^2}a_i, \quad (6)$$

$h$ is the number of non-zero terms in the sum $\sum_{i=1}^{m}\rho_i^+(x)$. Then, the pseudo-projection $\pi_M(x)$ of the point $x$ on the polytope $M$ is defined by the following equation:

$$\pi_M(x) = \lim_{k \to \infty}\varphi_M^{(k)}(x), \quad (7)$$

where

$$\varphi_M^{(k)}(x) = \underbrace{\varphi_M(\varphi_M(\ldots\varphi_M(x)\ldots))}_{k}. \quad (8)$$

```
Algorithm 1.
1:  input x_0
2:  k := 0
3:  x_{k+1} := x_{k+1} - φ_M(x_k)
4:  if ‖x_{k+1} - x_k‖ < ε go to 7
5:  k := k + 1
6:  go to 3
7:  output x̃ = x_{k+1}
8:  stop
```

Fig. 1. Algorithm computing pseudo-projection.

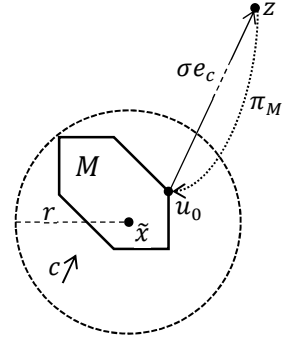

Fig. 2. Constructing the first approximation.

The iterative Algorithm 1 computing the pseudo-projection operation is shown in Fig. 1. The algorithm terminates when the distance between adjacent approximations becomes less than the small positive quantity $\varepsilon$ being a parameter of the algorithm. The proof of convergence of Algorithm 1 can be found in [30]. In [31], the authors proposed and investigated a scalable parallel implementation of Algorithm 1 in the form of operations on lists. By using the cost metric of the BSF parallel computation model, the authors showed that the scalability boundary[1] of given algorithm can be estimated as $O(\sqrt{n})$.

The *Target stage* calculates a sequence of points $\{u_0, u_1, \ldots, u_k, \ldots\}$ that satisfies the properties (2)-(4) and converges to the exact solution of the problem (1). To get the initial approximation $u_0$, we use the apex point $z$, which is defined as follows. Let $\tilde{x} \in M$ be the point obtained in the Quest stage. Let $S$ designate $n$-ball of radius $r$ and center $\tilde{x}$ satisfying the following property: $M \subset S$ (see Fig. 2). Let us fix some positive number $\sigma \in \mathbb{R}_{>0}$ such that $\sigma \gg r$. Let us define the unit vector

$$e_c = \frac{c}{\|c\|}, \quad (9)$$

where $c$ is the vector of the objective function coefficients of the problem (1). Then the *apex point* is calculated using the following equation:

$$z = \tilde{x} + \sigma e_c. \quad (10)$$

Let us set the initial approximation $u_0$ as follows:

$$u_0 = \pi_M(z). \quad (11)$$

In other words, $u_0$ is the pseudo-projection of the apex point $z$ on the polytope $M$. Using the proof scheme of the algorithm convergence in [30], it is easy to show, in this case, that the point $u_0$ will lie on the boundary $\Gamma_M$ of polytope $M$. This means that there is a number $i \in \{1, \ldots, m\}$ such that

$$u_0 \in H_i \cap \Gamma_M, \quad (12)$$

where $H_i$ is the hyperplane defined by the equation $\langle a_i, x \rangle = b_i$. Thus, condition (2) holds for the point $u_0$.

Now, let us assume that an approximation $u_k$ satisfying conditions (2) and (3) has already been found. To construct the next approximation $u_{k+1}$, let us calculate the intermediate point

$$v_k = u_k + \delta e_c, \quad (13)$$

that is obtained by adding to the point $u_k$ the unit vector $e_c$ multiplied by a small positive quantity $\delta \in \mathbb{R}_{>0}$ (see Fig. 3). Applying the pseudo-projection operation $\pi_M$ to the point $v_k$, we get the following intermediate point:

$$w_k = \pi_M(v_k). \quad (14)$$

Since the mapping $\varphi_M$ used to calculate the pseudo-projection by the equation (7) is a single-valued $M$-Fejerian mapping[2] [32], it follows from (12) and (13) that $w_k \in H_i \cap \Gamma_M$ for a sufficiently small $\delta$, i.e. $w_k$ will lie on the same face as $u_k$. If the value of objective function at the point $u_k$ is greater than or equal to the value of objective function at the point $w_k$, then the point $u_k$ is a solution of the problem (1). Let us assume that the opposite holds, i.e. $\langle c, w_k \rangle > \langle c, u_k \rangle$. Let $L_{u_k w_k}$ be the ray from the point $u_k$ passing through the point $w_k$:

$$L_{u_k w_k} = \{x \in \mathbb{R}^n \mid x = u_k + \eta(w_k - u_k), \eta \in \mathbb{R}_{\geq 0}\}.$$

Let us define the mapping $\gamma : \mathbb{R}^n \times \mathbb{R}^n \to \mathbb{R}^n$ as follows:

$$\gamma(u_k, w_k) = \arg\max_x \{\|x - u_k\| \mid x \in L_{u_k w_k} \cap M\}. \quad (15)$$

---

[1] The scalability boundary of a parallel algorithm for a cluster computing system is the maximum number of computing nodes up to which the speedup increases.

[2] A single-valued mapping $\varphi_M : \mathbb{R}^n \to \mathbb{R}^n$ is called *Fejerian with respect to the set* $M$, or briefly *M-Fejerian* if $\forall y \in M (\varphi(y) = y) \lor \forall x \notin M (\forall y \in M (\|\varphi(x) - y\| < \|x - y\|))$.

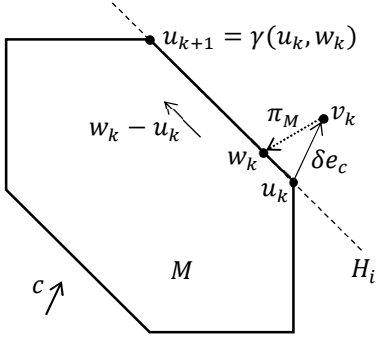

Fig. 3. Constructing the next approximation $u_{k+1}$.

In other words, the mapping $\gamma$ gives a point which lies on the ray $L_{u_k w_k}$, belongs to the polytope $M$, and is located as far from the point $u_k$ as possible. Let us use the given point as the next approximation:

$$u_{k+1} = \gamma(u_k, w_k). \qquad (16)$$

Obviously, the conditions (2) and (3) also hold for the point $u_{k+1}$.

The Algorithm 2 implementing the Target stage is shown in Fig. 4. Step 1 reads the point $\tilde{x} \in M$ obtained in the Quest stage. By using the vector $c$ consisting of the objective function coefficients of the problem (1), Step 2 computes the unit vector $e_c$ such that $e_c$ and $c$ are like vectors. Step 3 computes the apex point $z$ by using the equation (10) (see Fig. 2). By using Algorithm 1, Step 4 computes the initial approximation $u_0$ being the pseudo-projection of the apex point $z$ on the polytope $M$ (see Fig. 2). Step 5 assigns the zero value to the iteration counter $k$. Step 6 starts the iterative loop that calculates the intermediate point $v_k$ using the equation (13). By using the equation (14), Step 7 calculates the intermediate point $w_k$ that lies on the same face of the polytope $M$ as the current approximation $u_k$ does (see Fig. 3). Step 8 checks the stopping criterion. If the objective function value at the point $w_k$ does not exceed the objective function value at the point $u_k$, the iterative process ends, and $u_k$ is output as an approximate solution of the problem (1). Otherwise, the point $w_k$ sets the direction $(w_k - u_k)$ in which the value of the objective function will increase. Step 9 calculates the feasible point $u_{k+1}$ at which the maximum of the objective function in the direction $(w_k - u_k)$ is reached. Step 10 increases the iteration counter $k$ by one. In Step 11, we go to Step 6 starting the next iteration. The convergence of Algorithm 2 immediately follows from the closure and boundedness of the polytope $M$. Thus, condition (4) holds.

```
Algorithm 2.
1:    input x̃
2:    e_c := c/‖c‖
3:    z := x̃ + σe_c
4:    u_0 := π_M(z)
5:    k := 0
6:    v_k := u_k + δe_c
7:    w_k := π_M(v_k)
8:    if ⟨c,w_k⟩ ≤ ⟨c,u_k⟩ go to 12
9:    u_{k+1} = γ(u_k, w_k)
10:   k := k + 1
11:   go to 6
12:   output x̄ = u_k
13:   stop
```

Fig. 4. Target stage implementation.

The Algorithm 3 computing an approximate value of the mapping $\gamma$ used in step 9 of Algorithm 2 is shown in Fig. 5. This mapping is defined by the equation (15). Algorithm 3 computes a sequence of points $\{t_0, t_1, \ldots, t_j, \ldots\}$ such that:

$$t_0 = u_k, \qquad (17)$$
$$\lim_{j \to \infty} t_j = \gamma(u_k, w_k). \qquad (18)$$

Step 1 reads the given values of $u_k$ and $w_k$. Step 2 computes the unit vector $e_{u_k w_k}$ such that $e_{u_k w_k}$ and $(w_k - u_k)$ are like vectors. Step 3 assigns the constant $\mu \in \mathbb{R}_{>0}$ being a parameter of the algorithm to the variable $\tau_0$ being the initial shift quantity, and assigns the vector $u_k$ to the vector $t_0$ being the initial approximation. Step 4 assigns the zero value to the variables $j$ and $i$ being the iteration counters. Step 5 checks the stopping criterion. If the shift quantity $\tau_i$ is less than the small positive quantity $\varepsilon$ being a parameter of the algorithm then the control is passed to Step 13 outputting the value $u_{k+1} = t_j$ as a result and stopping the computations. Otherwise, Step 6 calculates the next approximation $t_{j+1}$ by shifting the point $t_j$ in the direction $(w_k - u_k)$ by the distance $\tau_i$. Step 7 checks if the point $t_{j+1}$ has gone outside the polytope $M$. If this is the case, we go to Step 10, which sets the new shift quantity $\tau_{i+1}$ to the value of $\frac{1}{2}\tau_i$. Step 11 increases the iteration counter $i$ by one and Step 12 passes the control to Step 6, which recalculates the point $t_{j+1}$ using the smaller shift quantity $\tau_i$. If the condition in Step 7 is false, i.e. the new approximation $t_{j+1}$ belongs to the polytope $M$, then Step 8 increases the counter $j$ by one and Step 9 passes the control to Step 5 (see Fig. 6). It is obvious that the sequence of points $\{t_0, t_1, \ldots, t_j, \ldots\}$ generated by Algorithm 3 satisfies the conditions (17) and (18).

```
Algorithm 3.
1:   input u_k, w_k
2:   e_{u_k w_k} := (w_k − u_k)/‖w_k − u_k‖
3:   τ_0 := μ;  t_0 := u_k
4:   j := 0;  i := 0
5:   if τ_i < ε go to 13
6:   t_{j+1} := t_j + τ_i e_{u_k w_k}
7:   if t_{j+1} ∉ M go to 10
8:   j := j + 1
9:   go to 6
10:  τ_{i+1} := τ_i/2
11:  i := i + 1
12:  go to 5
13:  output u_{k+1} = t_j; stop
```

Fig. 5. Computation of the mapping $\gamma$ defined by (15).

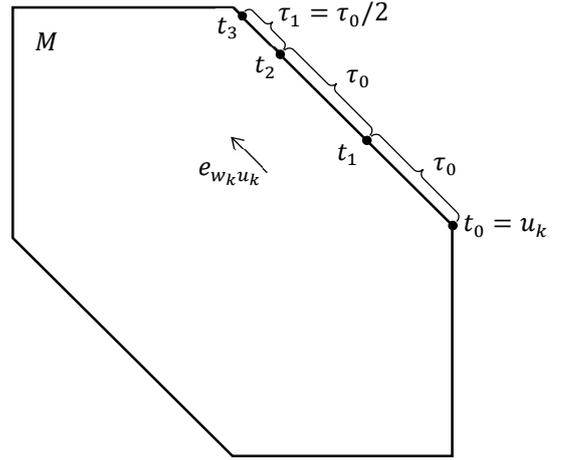

Fig. 6. Illustration of Algorithm 3 working principle.

## III. Software implementation and computational experiments

We performed a parallel implementation of Algorithm 2 in C++ using the BSF skeleton [33] based on the MPI parallel programming library. The parallel implementation scheme is similar to the ones used in [34] and [31]. The source codes are freely available at https://github.com/leonid-sokolinsky/NSLP-Quest. As a test problem, we used the scalable inequality system (19) of dimension $n$ from [28].

$$\begin{cases} x_0 & \leq 200 \\ \quad x_1 & \leq 200 \\ \quad \ddots & \ldots \ldots \\ \quad\quad x_{n-1} & \leq 200 \\ x_0 + x_1 \cdots + x_{n-1} & \leq 200(n-1)+100 \\ x_0 + x_1 \cdots + x_{n-1} & \geq 100 \\ x_0 & \geq 0 \\ \quad x_1 & \geq 0 \\ \quad \ddots & \ldots \ldots \\ \quad\quad x_{n-1} & \geq 0 \end{cases} \quad (19)$$

The number of inequalities in this system is $m = 2n + 2$. We used the following vector of the objective function coefficients:

$$c = (10n, 10(n-1), 10(n-2), \ldots, 10). \quad (20)$$

In this case, the exact solution of the problem (1) for any dimension $n \geq 2$ will be the point

$$\bar{x} = (200, 200, \ldots, 200, 100).$$

The computational experiments were conducted on the "Tornado SUSU" computing cluster [35], whose specifications are shown in Fig. 7. The results of the experiments are presented in Fig. 8. The computations were performed with the following dimensions: 5 000, 7 500 and 10 000. The number of inequalities was 10 002, 15 002 and 20 002, respectively.

Experiments have shown that the scalability boundary of the apex method depends significantly on the size of the problem. At $n = 5000$, the scalability boundary was approximately 55 processor nodes. For the problem of dimension $n = 7500$, this boundary increased to 80 nodes, and for the problem of dimension $n = 10000$, it was close to 100 nodes. Further increasing the problem size caused the compiler error: "insufficient memory". It should be noted that the computations were performed in double precision floating-point format occupying 64 bits in computer memory. An attempt to use single-precision floating-point format occupying 32 bits failed because Algorithm 1, which calculates the pseudo-projection, stopped to converge. The experiments also showed that the parameter $\sigma$ that determines the distance of the apex point $z$ from the polytope $M$ by the equation (10) (see Fig. 2) has little effect on the total time of solving the problem when this parameter has large values (more then 100000). For the LP problem with the constraints (19) and the objective function coefficients (20), the parameter $\sigma$ cannot be less than $200n$, since in this case the apex point $z$ will be inside the polytope $M$. If the apex point is not far enough away from the polytope, then its pseudo-projection can be an interior point of some polytope face. If the apex point is taken far enough away from the polytope (the value $\sigma = 20000n$ was used in the experiments), then the pseudo-projection always fell into one of the polytope vertices. Also, we would like to note that all computed points in the sequence $u_0, u_1, u_2, \ldots$ were vertices of the polytope.

The computational experiment showed that more than 99% of the time spent for solving the LP problem by the apex method was taken by the calculation of pseudo-projections (Step 4 of Algorithm 2). At that, the calculation of one approximation $u_k$ for a problem of dimension $n = 10000$ on 100 processor nodes took 44 minutes.

| Number of processor nodes | 480 |
|---|---|
| Processor | Intel Xeon X5680 (6 cores 3.33 GHz) |
| Processors per node | 2 |
| Memory per node | 24 GB DDR3 |
| Interconnect | InfniBand QDR (40 Gbit/s) |
| Operating system | Linux CentOS |

Fig. 7. Specifications of "Tornado SUSU" computing cluster.

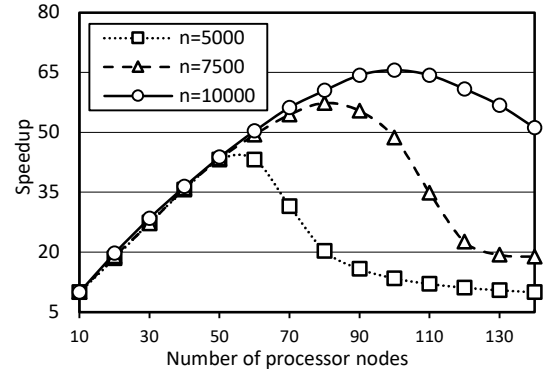

Fig. 8. Scalability of Algorithm 2.

## IV. Conclusion

The article proposes a new scalable iterative method called the "apex method" for solving the large linear programming problems. The apex-method uses the predictor-corrector framework and proceeds in two stages. The first stage called *Quest* searches for a feasible point of the linear programming problem. The second stage called *Target* calculates a sequence of points, which lie on the faces of polytope bounding the feasible region, and which converge to the exact solution of the LP problem. The presented method is implemented in C++ using the MPI parallel programming library. The computational experiments on solving a large linear programming problem on a cluster computing system are described. The conducted experiments showed that the apex method scales up well with increasing problem size. The strength of the apex method is its "self-correcting" in the presence of errors caused by loss of precision. The apex method can also potentially be used for solving non-stationary linear programming problems. The weakness of this method is the high computational complexity of constructing a pseudo-projection.

As future research directions, we suggest the following.

1. Perform a formal proof of the apex method convergence.
2. Execute testing of the apex method on random linear programming problems generated by a special algorithm.
3. Solve some linear programming problems from the Netlib-Lp repository [36], [37] using the apex method.
4. Use the apex method to generate a test dataset for the development and training of neural networks that can quickly solve large-scale linear programming problems in cooperation with a supercomputer.

## Acknowledgment

This research was partially supported by the Russian Foundation for Basic Research (project No. 20 07 00092 a), by the Ministry of Science and Higher Education of the Russian Federation (gov. order FENU-2020-0022) and by the Government of the Russian Federation according to Act 211 (contract No. 02.A03.21.0011).